\documentclass[12pt,a4paper]{article}
\usepackage{mathrsfs}
\usepackage{amsfonts}
\usepackage{cases}
\usepackage{times}
\usepackage{dsfont}
\usepackage{color}
\usepackage{graphics}
\newenvironment{keywords}{{\bf Keywords: }}{}

\begin{document}
\author{Hanbing Liu\footnote{Email address: liu.hanbing@uaic.ro}
\ \\
{\small Faculty of Mathematics, ``Alexandru Ioan Cuza'' University,
Carol I Bd., No.9-11, Iasi 700506, Romania}}

\title{Optimal control problems with state constraint governed by Navier-Stokes
equations\thanks {This work was supported by founds from the Marie
Curie ITN ``Controlled Systems'', No.213841/2008.}}
\date{}

\maketitle \noindent
\hrule
\begin{abstract} This work deals with the existence of optimal
solution and the maximum principle for optimal control problem
 governed by Navier-Stokes equations with state constraint in 3-D. Strong results in 2-D also are given.
\end{abstract}

\noindent
\ \\
\hrule

\bigskip
\section{Introduction}
In this paper, we shall study the optimal control problem \\
$$(\textbf{P})\ \ \ \ \textbf{Minimize}\ \
\frac{1}{2}\int_{0}^{T}\left(\int_{\Omega}|\mathscr{C}(y(t,x)-y^{0}(t,x))|^{2}\right)dxdt+\int_{0}^{T}h(u(t))dt;$$
subject to
\begin{equation}
\left\{\begin{array}{lllll}
     \frac{\partial y}{\partial t}-\nu \triangle y+(y\cdot\nabla)y+\nabla p=D_0u(t)+f_0(t),\ \ \ \ \ \mathrm{in}\ \Omega\times(0,T),\\
     y(0)=y_{0}\ \ \ \ \ \ \ \ \ \ \ \ \ \ \ \ \ \ \ \ \ \ \ \ \ \ \ \ \ \ \ \ \ \ \ \ \ \ \ \ \  \ \ \ \ \ \ \ \ \ \ \ \ \ \ \ \ \ \ \ \  \ \ \ \ \ \  \mathrm{in}\ \Omega,\\
     \nabla\cdot y=0\ \ \ \ \ \  \ \ \ \ \ \ \ \ \ \ \ \ \ \ \ \ \ \ \ \ \ \ \ \ \ \ \ \ \ \ \ \ \ \ \ \ \ \ \ \ \ \ \  \ \ \ \ \ \ \ \ \ \ \ \ \ \ \ \ \ \ \   \mathrm{in}\ \Omega\times(0,T),\\
     y=0\ \ \ \ \ \ \ \ \ \ \ \ \ \ \ \ \ \ \ \ \ \ \ \ \ \ \ \ \ \ \ \ \ \ \ \ \ \ \ \ \ \ \ \ \ \ \ \ \ \  \ \ \ \  \ \ \ \ \ \ \ \ \ \ \ \ \ \ \ \ \ \ \ \  \mathrm{on}\ \partial\Omega\times(0,T)
\end{array}\right.
\end{equation}
\begin{equation}
y(t)\in K,\ \ \ \ \ \ \ \ \forall t\in (0,T),
\end{equation}
where $K$ is a closed convex subset in
\begin{equation}
H=\{y;y\in (L^{2}(\Omega))^{N},\nabla\cdot y=0,y\cdot\mathbf{n}=0\
\mathrm{on}\ \partial\Omega\}.
\end{equation}
Here $\Omega$ is a bounded and open subset of $\mathbb{R}^N$ with
smooth boundary $\partial\Omega$, $T>0$ is a given constant, $\nu>0$
is the viscosity constant, $f_0\in L^2(0,T;(L^2(\Omega))^N)$ is a
source field, $y(x,t)$ is the velocity vector, $p$ stands for the
pressure, $D_0\in L(U;(L^2(\Omega))^N)$,  and $u\in L^2(0,T;U)$,
where $U$ is a Hilbert space.

The function $h:U\rightarrow (-\infty,+\infty]$ is convex and lower
semicontinuous, $y^0\in L^2(0,T;H)$, and $\mathscr{C}\in L(V,H)$,
where $V=((H_0^1(\Omega)))^N\cap H$. Two cases of physical interest
are covered by the
cost functional of this form,\\
(a)$\mathscr{C}=d_1I$: the objective is to minimize the energy of
control and the difference between state function and object
function; \\
(b)$\mathscr{C}=d_2\nabla\times$: the objective in this case is to
minimize the energy of control and regularize the smoothness
of the state function. \\
Physically, the cost functional in case (a) means a regulation of
turbulent kinetic energy, while in case (b) it means a regulation of
the square of vorticity. (See [6,7,8] for a discussion on this
control problem.)

Let us introduce some functional spaces and some
operators to represent the Navier-Stokes equation (1.1) as infinite
dimensional differential equations.

Denote by the symbol $\parallel \cdot
\parallel$ the norm of the space $V$, which is defined by
$$\parallel y\parallel^2=\sum_{i=1}^N\int_\Omega|\nabla y_i|^2dx,\ \ $$
and by the symbol $|\cdot|$ the norm of $\mathbb{R}^N$ and
$(L^{2}(\Omega))^{N}$. We endow the space $H$ with the norm of
$(L^{2}(\Omega))^{N}$, and denote by $\langle \cdot,\cdot\rangle$
the scalar product of $H$, $\langle \cdot,\cdot\rangle_{(V,V')}$ the
paring between $V$ and its dual $V'$ with the norm $\parallel
\cdot\parallel_{V'}$. Let $A\in L(V,V')$ and $b:V\times V\times
V\rightarrow\mathbb{R}$ be defined by:
$$\langle Ay,z\rangle=\sum_{i=1}^N\int_\Omega\nabla y_i\cdot\nabla z_idx,\ \forall y,z\in V$$
and
$$b(y,z,w)=\sum_{i=1}^N\int_\Omega y_iD_iz_jw_jdx,\ \forall y,z,w\in V$$
respectively, where $D_i=\frac{\partial}{\partial x_i}$.
$D(A)=(H^2(\Omega))^N\cap V$. We define $B:V\rightarrow V'$ by
$$\langle B(y),w\rangle=b(y,y,w),\ \forall y,w\in V$$
Let $f(t)=Pf_0(t)$ and $D\in L(U,H)$ be given by $D=PD_0$, where
$P:(L^{2}(\Omega))^{N}\rightarrow H$ is the projection on $H$. Then
we may rewrite the optimal control problem  $(P)$ as:
$$\textbf{(P)}\ \ \ \ \ \ \ \ \ \ \ \textbf{Min} \ \frac{1}{2}\int_{0}^{T}|\mathscr{C}(y(t)-y^{0}(t))|^{2}+\int_{0}^{T}h(u(t))dt;$$
subject to
\begin{equation}
\left\{\begin{array}{lllll}
     y'(t)+\nu Ay(t)+By(t)=Du(t)+f(t),\\
     y(0)=y_{0},
\end{array}\right.
\end{equation}
with
\begin{equation}
    y(t)\in K\ \ \ \ \ \ \ \ \ \ \forall t\in[0,T]
\end{equation}

Since $f, Du\in L^2(0,T;H), y_0\in V$, equation (1.4) has a unique
solution $y\in W^{1,2}(0,T;H)\cap L^2(0,T;D(A))$ when $N=2$ while in
the case $N=3$, for each $u\in L^2(0,T;U)$, there exists $0<T(u)\leq
T$ such that (1.4) has a unique solution $y(\cdot;u)\in
W^{1,2}(0,T^*;H)\cap L^2(0,T^*;D(A))$ for all $T^*<T(u)$. Here
$T(u)$ is given by
\begin{equation}
    T(u)=\frac{\nu}{3C_0^3[\|y_0\|^2+(\frac{1}{\nu})\|f+Du\|^2_{L^2(0,T;H)}]^3}
\end{equation}
where $C_0$ is a positive constant independent of $y_0,u$ and $\nu$
(see [3], p.261, Th.5.10).  In order to formulate the optimal
control problem governed by such system in terms of strong state
$y(\cdot;u)$, we observe from (1.6) that for each $L>0$, there
exists $T(L)>0$, such that for any $T^*<T(L)$ and any $u\in
L^2(0,T;U)$ with $\|Du\|_{L^2(0,T;H)}\leq L$, equation (1.4) has a
unique solution $y(\cdot;u)\in W^{1,2}(0,T^*;H)\cap
L^2(0,T^*;D(A))$. Therefore, the optimal control problem is
well-posed in the sense of strong solutions if we consider the
admissible control set as a bounded subset of $L^2(0,T;U)$. Another
way to formulate the control problem is in the framework of weak
solutions to equation (1.4), that is $y\in\mathscr{Y}_w=
L^2(0,T;V')\cap C_w(0,T;H)\cap W^{1,1}(0,T;V')$ , satisfying (see
[3],p.265,Th.5.12), for each $\Psi\in V$
\begin{equation}
\left\{\begin{array}{lllll}
     \frac{d}{dt}\langle y(t),\Psi\rangle_{(V',V)}+\nu a(y,\Psi)+b(y,y,\Psi)=\langle Du+f,\Psi\rangle_{(V',V)},\ a.e.\ t\in(0,T)\\
     y(0)=y_{0}.             \  \ \ \  \ \ \ \ \ \ \ \ \ \ \ \ \ \ \
\end{array}\right.
\end{equation}
where $C_w(0,T;H)$ is the space of weak continuous functions
$y:[0,T]\rightarrow H$. It is known that there exists at least a
weak solution to equation (1.4) for each $u\in L^2(0,T;U)$ (see
[3],p.265,Th5.12). We shall denote
$\mathscr{P}_w=\{(y,u)\in\mathscr{Y}_w\times L^2(0,T;U);(y,u)$
solution to (2.4)$, y(t)\in K,\forall t\in[0,T]\}$.

The main results of this work are about the existence of optimal
solution and maximum principle for problem (P) in 3-D. In [5,7,11],
some existence results are given for optimal control problems
governed by Navier-Stokes equations, wherein the admissible state
functions are considered as the strong solutions to Navier-Stokes
equations while in the present work, we give the existence result in
the framework of weak solutions to equation (1.4). In [5,10,11],
some Pontryagin's maximum principle type results are derived for
optimal control problems governed by Navier-Stokes equations. The
main differences between the present work and works mentioned above
are as follows. In this paper, we shall give the maximum principle
for problem (P) with state constraint of pointwise type, i.e. (1.5),
and it is not studied in [5,10,11], wherein the types of state
constraint involved include type of integral, type of two point
boundary and periodic type. We shall not only consider the stat
constraint set $K$ as a closed convex subset of $H$, but also derive
the maximum principle when $K$ is a closed convex subset of $V$.
Since the state constraint in the second case is stronger, the proof
is more precise, and the corresponding result is weaker, but
physically, it can be applied in some important examples in fluid
mechanics which will be given in section 4. This is also one
advantage of the results derived in this paper over those in the
mentioned works.

The outline of this paper is as follows. In section 2, we give and
prove the existence of the optimal pair for problem (P) by
considering the weak admissible pair set $\mathscr{P}_w$. In section
3, we shall formulate the optimal control problem in terms of strong
state function, which is different from that in section 2, and we
get the first order necessary conditions for problem (P) with state
constraint in two different cases mentioned above respectively. In
section 4, we give some examples of state constraint covered by the
two cases.

The following hypothesis will be in effected throughout this paper: \\
(i) $K\subset H$ is a closed convex subset with nonempty
interior; \\
(ii) $\mathscr{C}\in L(V;H)$, $D\in L(U;H)$, $y^0\in
L^2(0,T;H\cap D(\mathscr{C}^*\mathscr{C})),\ f\in L^2(0,T;H), \ y_0\in V$;\\
(iii) $h:U\rightarrow (-\infty,+\infty]$ is a convex lower
semicontinuous function. Moreover, there exist $\alpha>0$ and
$C\in\mathbb{R}$ such that
\begin{equation}
h(u)\geq \alpha|u|_U^2+C,\ \forall u\in U.
\end{equation}
When we study problem (P) in the case that $K$ is a closed convex
subset of $V$, we need assumption (ii') which is assumption (ii)
together with the assumption $D\in L(U;V)$.

We recall some properties of $b(y,z,w)$ here (see details in [3,9]):\\
$b(y,z,w)=-b(y,w,z)$, and there exists a positive constant $C$, such
that
$$|b(y,z,w)|\leq C\|y\|_{m_1}\|z\|_{m_2+1}\|w\|_{m_3}$$
where $m_1,m_2,m_3$ are positive numbers, satisfying:
$$\left\{\begin{array}{lllll} m_1+m_2+m_3\geq \frac{N}{2}, \
\mathrm{if}\ m_i\neq
\frac{N}{2},\forall i\in \{1,2,3\} \\
m_1+m_2+m_3> \frac{N}{2}, \ \mathrm{if}\ \exists i\in
\{1,2,3\},m_i=\frac{N}{2}
\end{array}\right.$$
We note also the interpolation inequality:
$$\|y\|_{m}\leq C\|y\|_{l}^{1-\alpha}\|y\|_{l+1}^{\alpha}$$
where $\alpha=m-l\in (0,1)$. Here $\|\cdot\|_{m}$ denotes the norm
of the Sobolev space $H^m(\Omega).$

We give some definition which will be used throughout this paper.\\
{\bf Definition 1.} Given a Banach space $E$ and its dual space
$E'$, we denote by $BV(0,T;E')$ the space of all functionals
$y:[0,T]\rightarrow E'$ with bounded variation. For each $\omega\in
BV(0,T;E')$, we define the continuous functional $\mu_\omega$ on
$C([0,T];E)$ by
\begin{equation}
\mu_\omega(z)=\int_0^T(z(t),d\omega(t))_{(E,E')}, \ \forall z\in
C([0,T];E)
\end{equation}
Here $(\cdot,\cdot)_{(E,E')}$ denotes the dual product between $E$
and $E'$, and the integral takes in the Riemann-Steiljes sense. The
measure $\mu_\omega$ will be denoted by $d\omega$, and if we assume
the space $E'$ is reflexive, then we have the Lebesgue decomposition
\begin{equation}
d\omega(t)=\omega_adt+d\omega_s(t)
\end{equation}
where $\omega_a\in L^1(0,T;E')$, $\omega_adt$ is the absolutely
continuous part of measure $d\omega$, and the functional
$\omega_s\in BV(0,T;E')$ is the singular part of $\omega$. In other
words, there exists a closed subset $\Theta\in[0,T]$ with the
Lebesgue measure zero such that $d\omega_s=0,\mathrm{on}\
[0,T]\setminus\Theta$ (see [2], p.51-p.57).

Denote by $M(0,T;E')$ the dual space of $C([0,T];E)$, i.e. the space
of all bounded $E'$-valued measures on [0,T], and notice that
$\mu_\omega\in M(0,T;E')$. We denote $\mathcal {K}=\{y\in
C([0,T];E);$ $ y(t)\in K, \forall t\in [0,T]\},$ and define the
normal cone to $\mathcal {K}$ at $y$ by
\begin{equation}
\mathcal {N}_\mathcal {K}(y)=\{\mu\in M(0,T;E');\mu(y-x)\geq
0,\forall x\in \mathcal {K}\}
\end{equation}

\section{Existence results}
 \setcounter{equation}{0} By admissible pair we mean $(y,u)\in \mathscr{P}_w$, which satisfies equation (1.4) in the weak sense, i.e.
 (1.7). An optimal pair is an admissible pair which minimizes $(P)$.\\
{\bf Theorem 1.} \emph{The optimal control problem (P) has at least
one optimal pair $(\hat{y},\hat{u})$. In 2-D, $\hat{y}$ is
strong solution to equation (1.4).} \\
Proof: When N=3, we denote
$$F(y,u)=\frac{1}{2}\int_{0}^{T}|\mathscr{C}(y(t)-y^{0}(t))|^{2}dt+\int_{0}^{T}h(u(t))dt$$
$$d_1=\inf\{\frac{1}{2}\int_{0}^{T}|\mathscr{C}(y(t)-y^{0}(t))|^{2}+\int_{0}^{T}h(u(t))dt;(y,u)\in\mathscr{P}_w\}.$$ Then there exist $(y_n,u_n)\in
\mathscr{P}_w$, such that
\begin{equation}
 d_1\leq F(u_n,y_n)\leq d_1+\frac{1}{n}.
\end{equation}
By (1.8) and (2.1), it follows that $\{u_n\}$ is bounded in
$L^2(0,T;U)$. Hence, there exists at least a subsequence which again
denoted by $\{u_n\}$, such that
\begin{equation}
u_n\rightarrow \hat{u} \ \mathrm{weakly\  in}\  L^2(0,T;U).
\end{equation}
Multiplying equation
\begin{equation}
 \noindent\left\{\begin{array}{lllll}
    y_n'(t)+\nu Ay_n(t)+By_n(t)=Du_n(t)+f(t),\\
    y_n(0)=y_{0}
\end{array}\right.
\end{equation}
by $y_n$, integrating on $(0,t)$, we get that
$$|y_{n}(t)|^{2}+\nu\int_{0}^{t}\|y_{n}\|^2ds\leq
C_1+C_2\int_{0}^{t}|y_{n}(s)|^2ds, $$  and it follows by Gronwall's
inequality that
\begin{equation}
|y_{n}(s)|^{2}+\nu\int_{0}^{T}\|y_{n}\|^2dt\leq C.
\end{equation}
This yields that
\begin{equation}
 \begin{array}{lllll}
  y_n\rightarrow\hat{y} \ \ \mathrm{ weak^*\ in}\ L^\infty(0,T;H),\ \mathrm{weakly\ in}\
  L^2(0,T;V),\\
  Ay_n\rightarrow A\hat{y}\ \mathrm{weakly\ in}\  L^2(0,T;V').
\end{array}
\end{equation}
By the properties of the trilinear function b, we have that
$$|\langle By_n,w\rangle_{(V',V)}|\leq C|y_n|^{\frac{1}{2}}\|y_n\|^{\frac{3}{2}}\|w\|,$$
and it follows that
\begin{equation}
\int_0^T|By_n|_{V'}^{\frac{4}{3}}dt\leq C\int_0^T\|y_n\|^2dt\leq C.
\end{equation}
Hence,
\begin{equation}
\int_0^T|\frac{dy_n}{dt}|_{V'}^{\frac{4}{3}}dt\leq C.
\end{equation}
Finally, we obtain by (2.6) and (2.7) that
\begin{eqnarray}
  \frac{dy_n}{dt}\rightarrow \frac{d\hat{y}}{dt}\ \mathrm{weakly\
  in}\
  L^{\frac{4}{3}}(0,T;V')\\
  By_n\rightarrow \eta\ \mathrm{weakly\ in}\ L^{\frac{4}{3}}(0,T;V')
\end{eqnarray}
To show that $(\hat{y},\hat{u})$ satisfies equation (1.4), it
remains to show that $\eta(t)=B\hat{y}(t),\ a.e.\ \mathrm{ in }\
(0,T).$ By (2.4), (2.7) and Aubin's compactness theorem(See[3],
p.26, Th.1.20), we obtain that
\begin{equation}
y_n\rightarrow \hat{y}\ \ \mathrm{strongly}\ \mathrm{in}\
L^2(0,T;H),
\end{equation}
and it follows that
\begin{equation}
 \begin{array}{lllll}
  \int_0^T|\langle
  By_n-B\hat{y},\psi\rangle_{(V',V)}|\leq\int_0^T\left(|b(y_n-\hat{y},y_n,\psi)|+|\hat{y},b(y_n-\hat{y},\psi)|\right)dt\\
  \rightarrow 0\ \ \ \mathrm{as}\ n\rightarrow +\infty,\ \forall\psi\in
  L^2(0,T;\mathscr{V}),
\end{array}
\end{equation}
where $\mathscr{V}=\{\psi\in C_0^\infty(\Omega); div\psi=0\}$.
Hence, $\eta(t)=B\hat{y}(t),\ a.e.\ \mathrm{ in }\ (0,T).$ Since $h$
is convex and lower semicontinuous, we obtain that
\begin{equation}
d_1\leq F(\hat{y},\hat{u})\leq \liminf_{n\rightarrow+\infty}
F(u_n,y_n)\leq d_1
\end{equation}
We also have that for each $t\in[0,T], \exists t_n\in(0,T),$ such
that $\hat{y}(t_n)\in K$, and
\begin{equation}
\hat{y}(t_n)\rightarrow \hat{y}(t)\ \ \ \mathrm{ weakly\  in} \ H.
\end{equation}
Since $K$ is a closed convex subset of $H$, it's weakly closed, and
this yields that $\hat{y}(t)\in K, \forall t\in [0,T]$. Hence,
$(\hat{y},\hat{u})$ is an optimal pair for problem $(P)$. $\sharp$\\
\ \ \\
{\bf Remark 1:} As we stated in Section 1, when $N=3$, if we assume
that the admissible control set is a bounded subset of $L^2(0,T;U)$,
then we can consider the strong solution in a local time interval
$(0,T^*)$. By the similar method applied in the proof of Theorem 1,
we can get the existence result, and the optimal state function
$\hat{y}\in W^{1,2}(0,T^*;H)\cap L^2(0,T^*;D(A))$. Moreover, the
same result follows when the state constraint set $K$ is a closed
convex subset of $V$.

\section{The maximum principle}\setcounter{equation}{0}
To get the maximum principle, we need to consider the strong
solution of the Navier-Stokes equations. As we mentioned in Section
2, when $N=3$, we need to consider the problem of such case with
bounded admissible control set $\mathscr{U}_{ad}=\{u\in
L^2(0,T;U);\|Du\|_{L^2(0,T;H)}\leq L\}$, and then we can consider
the strong solution to Navier-Stokes equation in $(0,T^*)$, where
$0<T^*=T(L+\delta)<T(L)$. Here $\delta>0$ is a fixed constant, and
$T(L)$ is given by (2.3), i.e.
\begin{equation}
    T(L)=\frac{\nu}{3C_0^3[\|y_0\|^2+(\frac{2}{\nu})(\|f\|_{L^2(0,T;H)}^2+L^2)]^3}
\end{equation}
Denote $\mathscr{D}(h)=\{u\in L^2(0,T;U);\int_0^Th(u)dt<+\infty\}$.
When $N=3$, we shall assume that
\begin{equation}
\mathscr{D}(h)\subset \mathscr{U}_{ad}.
\end{equation}
With this assumption, we can consider the strong solution in
$[0,T^*]$ in 3-D without control constraint which is included in the
definition of the function $h$ inexplicitly. We give an example here
to show that this assumption can be easily fullfilled. Let $h(u)$ be
$$h(u)=\left\{\begin{array}{lllll}
0,\ \ \ \ \ \ \mathrm{if}\ \|u\|_U\leq r;\\
+\infty,\ \mathrm{if}\ \|u\|_U>r
\end{array}\right. $$
We see that if $\int_0^Th(u(t))dt<+\infty$, then $\|u(t)\|_U\leq r,
a.e.\ \ \mathrm{in} \ [0,T]$, so $\|u\|_{L^2(0,T;U)}\leq C(r)$.
Since in 2-D, the strong solution to equation (1.4) exists on
arbitrary time interval $(0,T)$, such assumption is unnecessary. We
still denote the interval $[0,T^*]$ where assumption (3.2) holds by
$[0,T]$.

We need also the following assumption:\\
(iv) There are $(\tilde{z}, \tilde{u})\in C(0,T;H)\times L^2(0,T;U)$
solution to equation
\begin{equation}
\left\{\begin{array}{lllll}
    \tilde{z}'(t)+\nu A\tilde{z}(t)+(B'(y^{*}(t)))\tilde{z}(t)=B(y^{*}(t))+D\tilde{u}(t)+f(t),\\
    \tilde{z}(0)=y_{0}
\end{array}\right.
\end{equation}
such that $\tilde{z}(t)\in$ int$K,$ for $t$ in a dense subset of
$[0,T]$.

Here $y^{*}(t)$ is the optimal state function for the optimal
control problem $(P)$. Inasmuch as
$$B(y^*)\in L^2(0,T;H)),\ |(B'(y^*)z,z)|\leq \frac{\nu}{4}\|z\|^2+C_\nu|z|^2,$$
we know that equation (3.3) has a solution $\tilde{z}\in
W^{1,2}(0,T;H)\cap L^2(0,T;D(A))$.\\
\ \\
{\bf Theorem 2.}\emph{ Let $(y^*(t),u^*(t))$ be the optimal pair for
the optimal control problem $(P)$. Then under assumptions
(i)$\sim$(iv), there are $p(t)\in L^{\infty}(0,T;H)$ and
$\omega(t)\in BV(0,T;H)$, such that:
\begin{equation}
D^*p(t)\in\partial h(u^*(t)) \ \ a.e. [0,T],
\end{equation}
\begin{equation}
p(t)=-\int_t^TU(s,t)(\mathscr{C}^*\mathscr{C}(y^{*}(t)-y^{0}(t)))ds-\int_t^TU(s,t)d\omega(s),
\end{equation}
and
\begin{equation}
\int_0^T\langle d\omega(t), y^*(t)-x(t)\rangle\geq 0,\forall x\in
\mathcal {K}
\end{equation}}
Here $D^{*}$, $\mathscr{C}^*$, $B'(y^{*}(t))^{*}$, are the adjoint
operators of $D$, $\mathscr{C}$ and $B'(y^{*}(t))$ respectively,
where $B'(y)$ is the operator defined by
$$\langle B'(y)z,
w\rangle=b(y,z,w)+b(z,y,w),\ \forall \ z,w\in V$$ We recognize in
(3.5) the mild form of the dual equation
\begin{equation}
\noindent\left\{\begin{array}{lllll}
     p'(t)=\nu Ap(t)+(B'(y^{*}(t))^{*})p(t)+\mathscr{C}^*\mathscr{C}(y^{*}(t)-y^{0}(t))+\mu_\omega(t),\ a.e.\ \mathrm{in}\ (0,T)\\
     p(T)=0
\end{array}\right.
\end{equation}
\ \\
Theorem 3 below is the analogue of Theorem 2 under the weaker
assumption : \\
(v) $K$ is a closed convex subset of $V$, and there
are $(\tilde{z}, \tilde{u})\in C(0,T;H)\times L^2(0,T;U)$ solution
to equation (3.3), such that $\tilde{z}(t)\in$ int$_VK,$ for $t$ in
a dense subset of $[0,T]$.

Here int$_VK$ is the
interior of $K$ with respect to topology of $V$.\\
{\bf Theorem 3.}\emph{ Let $(y^*(t),u^*(t))$ be the solution for
optimal control problem $(P)$. Then under assumptions
(ii'),(iii),(v), there are $p(t)\in L^{\infty}(0,T;V'), \omega(t)\in
BV(0,T;V')$, such that (3.4) and (3.5) hold, and (3.6) holds in the
sense of
\begin{equation}
\int_0^T\langle d\omega(t), y^*(t)-x(t)\rangle_{(V',V)}\geq
0,\forall x\in \mathcal {K}
\end{equation}}
\ \\
We define first the approximating cost functions
\begin{equation}F_{\lambda}(y,u)=\frac{1}
{2}\int_{0}^{T}\left(|\mathscr{C}(y-y^{0})|^{2}+|u-u^{*}|^{2}_U\right)dt+\int_{0}^{T}\left(\varphi_{\lambda}(y)+h_{\lambda}(u)\right)dt.
\end{equation}
where $h_{\lambda}$ and $\varphi_{\lambda}$ are the regularizations
of $h$ and $\varphi$ respectively, that is
\begin{equation}
h_{\lambda}(u)=\inf\{\frac{\|u-v\|_U^{2}}{2\lambda}+h(v); v\in U\},\
\ \varphi_{\lambda}(y)=\inf\{\frac{|y-x|^{2}}{2\lambda}+\varphi(x);
x\in H\}
\end{equation}
Here $\varphi$ is the characteristic function of $K$, which is
defined by
$$
\varphi(x)=\left\{\begin{array}{lllll}
+\infty\ \ \ \ \ \ \ \ \ \ \ \ \mathrm{if \ }x\in H\setminus K\\
0 \ \ \ \ \ \ \ \  \ \ \ \ \ \ \ \ \ \mathrm{if}\ \ x\in K.\\
\end{array}
\right. $$ The function $\varphi_{\lambda}$ is convex, continuous,
Gateaux differentiable, and
$\partial\varphi_{\lambda}=\nabla\varphi_{\lambda}=(\partial\varphi)_{\lambda}$,
which is single-valued (see details in [3], p.48, Th.2.9). Denote
$$\mathscr{P}=\{(y,u)\in C([0,T];H)\times
L^2(0,T;U);(y(t),u(t)) \mathrm{\ satisfy\  equation }\  (1.4)\}.$$
We prove first\\
{\bf Lemma 1.} \emph{There exists at least one optimal pair
$(y_\lambda,u_\lambda)$ for the optimal control problem:
$$(P_\lambda)\ \ \ \ \ \ \ \ \ \ \ \ \ \mathrm{Min}\{F_{\lambda}(y,u); (y,u)\in\mathscr{P}\}. $$
and $(y_\lambda,u_\lambda)\rightarrow (y^*,u^*)$ strongly in
$C([0,T];H)\cap L^2(0,T;V)\times L^2(0,T;U)$. Moreover, $\{y_\lambda\}$ is bounded in $C([0,T];V)\cap L^2(0,T;D(A))$.} \\
{\bf Proof:} The existence of the optimal pair follows by Theorem 1
and the arguments in remark 1. We shall show the convergence of the
optimal pair $(y_\lambda,u_\lambda)$ in 3-D, and it's easy to prove
that the same results hold in 2-D by applying the similar method.
Since
\begin{equation}
d_\lambda=F_\lambda(y_\lambda, u_\lambda)\leq F_\lambda(y^*,
u^*)\leq F(y^*, u^*)=d,
\end{equation}
we have that $\int_0^Th_\lambda(u_\lambda)dt\leq C, \forall
\lambda>0$. Since
$$h(J^h_\lambda(u_\lambda))\leq h_\lambda(u_\lambda)\leq h(u_\lambda), $$
where $J^h_\lambda=(1+\lambda\partial h)^{-1}$, we know that
$\int_0^Th(J^h_\lambda(u_\lambda))\leq C, \forall \lambda>0$, and by
assumption (3.2), we obtain that
$\|D(J^h_\lambda(u_\lambda))\|_{L^2(0,T;H)}\leq L, \forall
\lambda>0$. Since
$$h_\lambda(u_\lambda)=\frac{\lambda}{2}\|\partial h_\lambda(u_\lambda)\|_U^2
+h(J^h_\lambda(u_\lambda))\geq \frac{\lambda}{2}|\partial
h_\lambda(u_\lambda)|^2+\alpha\|J^h_\lambda(u_\lambda)\|_U^2+C,$$ we
have that $\{\lambda \|\partial h_\lambda(u_\lambda)\|^2\}$ is
bounded in $L^1(0,T)$, and it follows that,
$$\int_0^T\|u_\lambda-J^h_\lambda(u_\lambda)\|_U^2dt=\lambda \int_0^T\lambda|\ \partial h_\lambda(u_\lambda)\|^2_Udt\rightarrow 0 \ \mathrm{as}\  \lambda\rightarrow 0.$$
This implies that $\exists\lambda_0>0$, such that
$\forall\lambda>\lambda_0 $
\begin{equation}
 \|Du_\lambda\|_{L^2(0,T;H)}\leq L+\frac{\delta}{2}.
\end{equation}
Multiplying equation
\begin{equation}
\left\{\begin{array}{lllll}
  y'_\lambda(t)+\nu Ay_\lambda(t)+By_\lambda(t)=Du_\lambda(t)+f(t)\\
  y_\lambda(0)=y_0
\end{array}\right.
\end{equation}
by $y_\lambda$, integrating on $(0,t)$, it follows by Gronwall's
inequality that
\begin{equation}
 |y_\lambda(t)|^{2}+\int_{0}^{T}\|y_\lambda(t)\|^2dt\leq
C,\ \forall \lambda>0.
\end{equation}
Multiplying equation (3.13) by $Ay_\lambda$, integrating on $(0,t)$,
with the inequality
$$
|\langle By_\lambda, Ay_\lambda\rangle|\leq
C\|y_\lambda\|^{\frac{3}{2}}|Ay_\lambda|^{\frac{3}{2}},
$$
we obtain that
$$
\|y_\lambda(t)\|^2+\frac{\nu}{2} \int_{0}^{T}|Ay_\lambda(t)|^2dt\leq
C_0\left(\|y_0\|^2+\frac{1}{\nu}\int_{0}^{T}|f+Du_\lambda|^2dt+\frac{1}{\nu}\int_{0}^{t}\|y_\lambda(s)\|^6ds\right).
$$
Here $C_0$ is the same constant as that in (3.1). It follows by
Gronwall's inequality that
$$\|y_\lambda(t)\|^2\leq \phi(t)£¬$$
where
\begin{eqnarray*} \phi(t)=\left(\frac{\nu
\phi^3(0)}{\nu-3t\phi^3(0)}\right)^{\frac{1}{3}}, \forall
t\in\left(0,\frac{\nu}{3\phi^3(0)}\right)\\
\phi(0)=C_0(\|y_0\|^2+\frac{1}{\nu}\int_{0}^{T}|f+Du_\lambda|^2dt)
\end{eqnarray*}
 By (3.12) and the
definition of $T$, we have that
$$T<T_1=\frac{\nu}{3C_0^3\left(\|y_0\|^2+\frac{2}{\nu}(|f|^2_{L^2(0,T;H)}+(L+\frac{\delta}{2})^2)\right)^3}\leq\frac{\nu}{3\phi^3(0)}.$$
Hence
\begin{equation}
\|y_\lambda(t)\|^2+\int_{0}^{T}|Ay_\lambda(t)|^2dt\leq C(\delta), \
\forall \lambda>\lambda_0.
\end{equation}
We mention here that $C(\delta)$ is a constant dependent on
$\delta$, and we shall denote all the constants by $C$ in the
following without emphasis. By the properties of $b$, and (3.14),
(3.15), we get
\begin{equation}
\|By_\lambda\|_{L^{2}(0,T;H)}\leq C\|Ay_\lambda\|_{L^{2}(0,T;H)}\leq
C,\ \forall \lambda>\lambda_0.
\end{equation}
This yields that
\begin{equation}
\|(y_\lambda)'\|_{L^{2}(0,T;H)}\leq C,\ \forall \lambda>\lambda_0.
\end{equation}
Therefore, on a subsequence convergent to 0, again denoted by
$\lambda$, we have
\begin{equation}
\begin{array}{lllll}
y_{\lambda}(t)\rightarrow y_1(t)\ \mathrm{strongly\ in}\
C([0,T;H])\cap
L^2(0,T;V)\\
Ay_\lambda(t)\rightarrow Ay_1(t) \ \mathrm{weakly \ in}\
L^{2}(0,T;H) \\
(y_\lambda(t))'\rightarrow y_1'(t) \ \mathrm{weakly\  in}\
L^{2}(0,T;H) \\
u_\lambda(t)\rightarrow u_1(t) \ \mathrm{weakly \ in} \ L^{2}(0,T;U)
\end{array}
\end{equation}
Since
$$| By_\lambda-By_1|\leq C\left(|y_\lambda-y_1|^{\frac{1}{2}}|Ay_\lambda|+\|y_\lambda-y_1\|^{\frac{1}{2}}|Ay_\lambda-Ay_1|^{\frac{1}{2}}\right),$$
by (3.15) and (3.18), we also have
\begin{equation}
By_\lambda(t)\rightarrow By_1(t) \ \mathrm{strongly \ in}\
L^{2}(0,T;H)
\end{equation}
So $(y_1(t),u_1(t))$ is a solution to equation (2.1). Moreover,
since
$$\varphi_\lambda(y_\lambda)=\frac{\lambda}{2}|\partial\varphi_\lambda(y_\lambda)|^2
+\varphi(J_\lambda^\varphi(y_\lambda))\geq
\frac{\lambda}{2}|\partial\varphi_\lambda(y_\lambda)|^2,$$ we know
that $\{\lambda|\partial\varphi_\lambda(y_\lambda)|^2 \}$ is bounded
in $L^1(0,T)$, and since
$\partial\varphi_\lambda(y_\lambda)=\frac{1}{\lambda}(y_\lambda-J_\lambda^\varphi(y_\lambda))$,
where $J_\lambda^\varphi(y_\lambda)$ is defined by the inclusion
$J_\lambda^\varphi(y_\lambda)-y_\lambda+\lambda\partial\varphi(J_\lambda^\varphi(y_\lambda))\ni
0$, we have
$$\int_0^T|y_\lambda-J_\lambda^\varphi(y_\lambda)|dt
\leq \lambda
T\int_0^T\lambda|\partial\varphi_\lambda(y_\lambda)|^2dt\rightarrow
0\ as\ \lambda\rightarrow 0,$$ so
$y_\lambda-J_\lambda^\varphi(y_\lambda)\rightarrow 0\ a.e.\ (0,T). $
Since $J_\lambda^\varphi(y_\lambda)\in K,\ \forall t\in[0,T]$, we
have that $y_1(t)\in K.\ \forall t$$\in [0,T]$. By (3.18) we know
that
\begin{equation}
\lim_{\lambda\rightarrow0}J_\lambda^h(u_\lambda)=u_1\ \mathrm{weakly
\ in\ } L^2(0,T;U).
\end{equation}
Since the convex function $u\rightarrow \int_0^Th(u)dt$ is lower
semicontinuous, we obtain that
\begin{equation}
\liminf_{\lambda\rightarrow0}\int_0^Th_\lambda(u_\lambda)\geq
\int_0^Th(u_1)dt.
\end{equation}
Inasmuch as
$$\liminf _{\lambda\rightarrow 0}F_{\lambda}(y_\lambda,
u_\lambda)\leq \lim_{\lambda\rightarrow 0}F_{\lambda}(y^*, u^*)\leq
F(y^*, u^*),$$ it follows by (3.18) and (3.21) that
$$\frac{1}
{2}\int_{0}^{T}\left(|\mathscr{C}(y_1(t)-y^{0}(t))|^{2}+|u_1(t)-u^{*}(t)|_U^{2}\right)dt+\int_{0}^{T}h(u_1(t))dt\leq\liminf
_{\lambda\rightarrow 0}F_{\lambda}(y_\lambda, u_\lambda)$$
$$\leq\frac{1}{2}\int_{0}^{T}\left(|\mathscr{C}(y^*(t)-y^{0}(t))|^{2}\right) dt+\int_{0}^{T}h(u^*(t))dt$$
$$\leq\frac{1}{2}\int_{0}^{T}\left(|\mathscr{C}(y_1(t)-y^{0}(t))|^{2}\right)dt+\int_{0}^{T}h(u_1(t))dt.$$
This yields that $u_1=u^*,\ y_1=y^*$ and $u_\lambda(t)\rightarrow u^*(t)$ strongly in $L^2(0,T;U)$.\\
\ \\
{\bf Lemma 2.} \emph{Let $z_\lambda(t)$ be the solution to the
equation:
\begin{equation}
\noindent\left\{\begin{array}{lllll}
     z_{\lambda}'(t)+\nu Az_{\lambda}(t)+(B'(y_{\lambda}(t)))z_{\lambda}(t)=B(y_{\lambda}(t))+D\tilde{u}(t)+f(t),\\
     z_{\lambda}(0)=y_{0}.
\end{array}\right.
\end{equation}
Then $z_\lambda(t)\rightarrow \tilde{z}(t)$ strongly in
$C([0,T];H)\cap L^2(0,T;V)$, where $(\tilde{z}(t),\tilde{u}(t))$ is
defined in equation (3.3), and $y_{\lambda}(t)$ is the the optimal
solution in lemma 1.}\\
{\bf Proof:}\ Multiplying equation (3.22) by $z_\lambda(t)$,we get
$$\frac{1}{2}\frac{d}{dt}|z_\lambda(t)|^2+\nu \|z_\lambda(t)\|^2+b(z_\lambda(t),y_\lambda(t),z_\lambda(t))
=b(y_\lambda(t),y_\lambda(t),z_\lambda(t))+\langle
f+D\tilde{u},z_\lambda\rangle$$ Integrating on $(0,t)$, since
$$|b(z_\lambda,y_\lambda,z_\lambda)|\leq C\|z_\lambda\|^{\frac{3}{2}}|z_\lambda|^{\frac{1}{2}}\|y_\lambda\|,
|b(y_\lambda,y_\lambda,z_\lambda|\leq
C\|y_\lambda\|^{\frac{3}{2}}|Ay_\lambda|^{\frac{1}{2}}|z_\lambda|,$$
we obtain by Young's inequality that
$$|z_\lambda(t)|^{2}+\frac{\nu}{2} \int_0^t\|z_\lambda(s)\|^2ds\leq C_1\int_0^t|z_\lambda(s)|^{2}ds+C_2. $$
It follows by Gronwall's inequality that
\begin{equation}
|z_\lambda(t)|^{2}+\int_{0}^{T}\|z_\lambda(t)\|^{2}dt\leq C.
\end{equation}
Multiplying equation (3.22) by $Az_\lambda(t)$, integrating from $0$
to $t$, we get that
$$\frac{1}{2}\|z_\lambda(t)\|^{2}-\frac{1}{2}\|y_0\|^2+\nu \int_0^t|Az_\lambda(s)|^2ds$$
$$=\int_0^tb(y_\lambda(s),z_\lambda(s),Az_\lambda(s))+
b(z_\lambda(s),y_\lambda(s),Az_\lambda(s))$$
$$+b(y_\lambda(s),y_\lambda(s),Az_\lambda(s))ds+\langle
f(s)+D\tilde{u},z_\lambda(s)\rangle ds.$$ Since
$$|b(y_\lambda,z_\lambda,Az_\lambda)+
b(z_\lambda,y_\lambda,Az_\lambda)|\leq
C(\|y_\lambda\|\|z_\lambda\|^{\frac{1}{2}}|Az_\lambda|^{\frac{3}{2}}+\|z_\lambda\|\|y_\lambda\|^{\frac{1}{2}}|Ay_\lambda|^{\frac{1}{2}}|Az_\lambda|),$$
$$|b(y_\lambda,y_\lambda,Az_\lambda)|\leq C\|y_\lambda\|^{\frac{3}{2}}|Ay_\lambda|^{\frac{1}{2}}|Az_\lambda|,$$
we obtain that
$$\frac{1}{2}\|z_\lambda(t)\|^{2}+\frac{\nu}{2} \int_0^t|Az_\lambda(s)|^2ds\leq C_1\int_0^t\|z_\lambda(s)\|^{4}ds+C_2.$$
It follows by Gronwall's inequality and (3.23) that
\begin{equation}
\|z_\lambda(t)\|^{2}+\int_0^T|Az_\lambda(t)|^2dt\leq C
\end{equation}
Since
$$|\langle(B'(y_{\lambda}(t)))z_{\lambda}(t),w\rangle|\leq
C\left(|Az_\lambda(t)|^\frac{1}{2}+|Ay_\lambda(t)|^\frac{1}{2}\right)|w|$$
we obtain by (3.23) and (3.24) that
\begin{equation}
\int_0^T|(B'(y_{\lambda}(t)))z_{\lambda}(t)|^2dt\leq C
\end{equation}
By (3.23), (3.24) and (3.25), we get that
\begin{equation}
\int_0^T|z'_\lambda(t)|^2dt\leq C
\end{equation}
Hence
\begin{equation}
\begin{array}{lllll}z_{\lambda}(t)\rightarrow \bar{z}(t)\ \mathrm{strongly}\ \mathrm{in}\ C([0,T;H])\cap
L^2(0,T;V)\\
Az_{\lambda}(t)\rightarrow A\bar{z}(t)\ \mathrm{weakly}\ \mathrm{in}\ L^2(0,T;H)\\
z'_\lambda(t)\rightarrow \bar{z}'(t)\ \mathrm{weakly}\ \mathrm{in}\
L^2(0,T;H).
\end{array}
\end{equation}
Since
\begin{equation}
\begin{array}{lllll}
|\langle(B'(y_{\lambda}(t)))z_{\lambda}(t)-(B'(y^*(t)))\bar{z}(t),w\rangle|\\
\leq
C(|y_{\lambda}-y^*|^{\frac{1}{2}}|Az_{\lambda}|+\|z_\lambda-\bar{z}\|^{\frac{1}{2}}
|A(z_{\lambda}-\bar{z})|^{\frac{1}{2}}\\
+|z_\lambda-\bar{z}|^{\frac{1}{2}}|Ay_{\lambda}|\|z_\lambda-\bar{z}\|^{\frac{1}{2}}+
\|y_{\lambda}-y^*\|^{\frac{1}{2}}|A(y_{\lambda}-y^*)|^{\frac{1}{2}})|w|
\end{array}
\end{equation}
and $y_{\lambda}\rightarrow y^*$, $z_\lambda\rightarrow\bar{z}$
strongly in $L^2(0,T;V)\cap C([0,T;H])$, we have also
\begin{equation}
B'(y_{\lambda}(t)))z_{\lambda}(t)\rightarrow (B'(y^*(t)))\bar{z}(t)\
\mathrm{strongly\ in}\ L^2(0,T;H)
\end{equation}
With above inequalities, passing $\lambda$ to $0$ in equation
(3.22), we find that $\bar{z}(t)$ satisfies the equation (3.3), and
by the
uniqueness of the solution, $\bar{z}(t)=\tilde{z}(t)$.$\sharp$\\

We shall denote by $U(t,s)$ and $U_\lambda(s,t), 0\leq s\leq t\leq
T$ the evolution operators generated by $\nu A+(B'(y^*(t)))^*$ and
$\nu A+(B'(y_\lambda(t)))^*$ respectively, which are given by
$U(t,s)\xi=\psi(t)$ and $U_\lambda(s,t)x=\psi_\lambda(t)$ for $0\leq
s\leq t\leq T$, where $\psi(t)$ and $\psi_\lambda(t)$ are the
solutions to
\begin{equation}
\left\{\begin{array}{lllll}
    \psi'(t)=\nu A\psi(t)+(B'(y^*(t)))^*\psi(t),\\
    \psi(s)=\xi
\end{array}\right.
\end{equation}
and
\begin{equation}
\left\{\begin{array}{lllll}
    \psi_\lambda'(t)=\nu A\psi_\lambda(t)+(B'(y_\lambda(t)))^*\psi_\lambda(t),\\
    \psi_\lambda(s)=\xi
\end{array}\right.
\end{equation}
respectively. It is well known and easily seen that such evolution
operators exist. Denote by $U^*(t,s)$, $U_\lambda^*(s,t)$ the
respective adjoint operators of $U(t,s)$ and $U_\lambda(s,t)$, which
are generated by $\nu A+B'(y^*(t))$ and $\nu A+B'(y_\lambda(t))$
respectively. By the similar method applied in lemma 2, we can
obtain that
\begin{equation}
\left\{\begin{array}{lllll}
     \|U_\lambda^*(t,s)\|_{L(H,H)}\leq C\\
     U_\lambda^*(s,t)\xi\rightarrow U^*(s,t)\xi\ \ \mathrm{in} \ C([0,T];H),
     \forall \xi\in H
\end{array}\right.
\end{equation}
\ \\
\ \\
{\bf Proof of Theorem 2:}\\
{\bf step 1:}(first order necessary condition for approximate
problem) Since $(y_\lambda, u_\lambda)$ minimize the functional
$F_\lambda(y, u), $ we know that
$$\lim_{\rho\rightarrow 0} \frac{F_\lambda(u_\lambda+\rho u)-F_\lambda(u_\lambda)}{\rho}=0,\ \ \forall u\in U,$$
and this yields
\begin{equation}
\langle \mathscr{C}^*\mathscr{C}(y_\lambda-y^0),w_\lambda \rangle+
(\nabla h_\lambda(u_\lambda)+u_\lambda-u^*,u)_U+\langle
\partial\varphi_\lambda(y_\lambda),w_\lambda\rangle=0,
\end{equation}
where $w_\lambda=\lim_{\rho\rightarrow 0}
\frac{y_\lambda^\rho-y_\lambda}{\rho}, \
(y_\lambda^\rho,u_\lambda+\rho u)\in \mathscr{P}$, and
$w_\lambda(t)$ is the solution to the equation
\begin{equation}
w'_\lambda(t)+\nu Aw_\lambda(t)+B'(y_\lambda(t))w_\lambda(t)=Du ,\
w_\lambda(0)=0.
\end{equation}
Let $p_\lambda(t)$ be the solution to the backward dual equation
\begin{equation}
\noindent\left\{\begin{array}{lllll} p'_\lambda(t)=\nu
Ap_\lambda(t)+(B'(y_\lambda(t))^*)p_\lambda(t)+\mathscr{C}^*\mathscr{C}(y_\lambda(t)-y^0(t))+\partial\varphi_\lambda(y_\lambda(t))\\
p_\lambda(T)=0
\end{array}\right.
\end{equation}
By (3.33),(3.34) and (3.35), we get by calculation that
\underline{}$$\langle p'_\lambda(t),w_\lambda(t) \rangle+ \langle
-Ap_\lambda(t)-(B'(y_\lambda(t))^*)p_\lambda(t),w_\lambda(t)\rangle+\langle
u_\lambda-u^*,u\rangle_U=0.$$ Hence
$$\langle-D^*p_\lambda(t)+\nabla h_\lambda(u_\lambda)+u_\lambda-u^*,u\rangle_U=0,\ \ \forall u\in U$$
Finally, we obtain that
\begin{equation}
D^*p_\lambda(t)=\nabla h_\lambda(u_\lambda)+u_\lambda-u^*(t), \ a.e.
\ t\in [0,T]
\end{equation}
\ \\
{\bf step2:}\ (pass $\partial\varphi_{\lambda}(y_\lambda),
p_\lambda, \partial h_\lambda(u_\lambda)$ to limit) By assumption
(iv) and lemma 2, we know that $\exists \rho>0, \lambda_1>0\ s.t.\
z_\lambda(t)+\rho h\in K,\ $ for $t$ in a dense subset of [0,T],\
$\forall |h|=1,\forall \lambda>\lambda_1.$ For $\lambda$ fixed,
$z_\lambda(t)$ is continuous in $[0,T]$, so there exists a partition
of $[0,T]$,
$$0=t_1<t_2<\cdots<t_{N-1}<t_N=T$$
such that $|z_\lambda(t_i)-z_\lambda(t_{i-1})|<\frac{\rho}{2},\
z_\lambda(t_i)+\rho h\in K,\forall 1\leq i\leq N$. Moreover, since
$$\sum_{i=1}^{N}\int_{t_{i-1}}^{t_i}\langle\partial\varphi_{\lambda}(y_\lambda(t)), y_\lambda(t)-z_\lambda(t_i)-\rho h\rangle dt\geq
\sum_{i=1}^{N}\int_{t_{i-1}}^{t_i}\varphi_{\lambda}(y_\lambda(t))-\varphi_{\lambda}(z_\lambda(t_i)+\rho
h)dt\geq0,$$ we get that
$$\rho \int_{0}^{T}|\partial\varphi_{\lambda}(y_\lambda)|dt\leq
 \sum_{i=1}^{N}\int_{t_{i-1}}^{t_i}\langle\partial\varphi_{\lambda}(y_\lambda(t)), y_\lambda(t)-z_\lambda(t_i)\rangle dt$$
$$=\int_{0}^{T}\langle\partial\varphi_{\lambda}(y_\lambda(t)), y_\lambda(t)-z_\lambda(t)\rangle dt
+\sum_{i=1}^{N}\int_{t_{i-1}}^{t_i}\langle\partial\varphi_{\lambda}(y_\lambda(t)),
z_\lambda(t)-z_\lambda(t_i)\rangle dt,$$ and it follows that
$$\frac{\rho}{2}\int_{0}^{T}|\partial\varphi_{\lambda}(y_\lambda)|dt\leq\int_{0}^{T}\langle\partial\varphi_{\lambda}(y_\lambda(t)), y_\lambda(t)-z_\lambda(t)\rangle dt$$
$$=\int_{0}^{T}\langle p'_\lambda(t)-\nu
Ap_\lambda(t)-(B'(y_\lambda(t))^*)p_\lambda(t)-\mathscr{C}^*\mathscr{C}(y_\lambda(t)-y^0(t)),
y_\lambda(t)-z_\lambda(t)\rangle dt$$
$$ =\int_{0}^{T}\langle p_\lambda,-y'_\lambda-\nu Ay_\lambda-(B'(y_\lambda))y_\lambda
+z'_\lambda+\nu Az_\lambda+(B'(y_\lambda))z_\lambda\rangle+\langle
\mathscr{C}^*\mathscr{C}(y_\lambda-y^0),y_\lambda-z_\lambda\rangle
dt$$
$$=\int_{0}^{T}\langle p_\lambda,D\tilde{u}-Du_\lambda\rangle-\langle
\mathscr{C}^*\mathscr{C}(y_\lambda-y^0),y_\lambda-z_\lambda\rangle
dt$$
$$=\int_{0}^{T}\left(\langle \nabla h_\lambda(u_\lambda)+u_\lambda(t)-u^*,\tilde{u}(t)-u_\lambda(t)\rangle_U-\langle
\mathscr{C}^*\mathscr{C}(y_\lambda-y^0),y_\lambda-z_\lambda\rangle\right)dt$$
\begin{equation}
\leq \int_{0}^{T}[h_\lambda(u_\lambda)-h_\lambda(\tilde{u})]dt+C\leq
C.
\end{equation}
We set
$\omega_\lambda(t)=\int_0^t\partial\varphi_{\lambda}(y_\lambda(s))ds,\
t\in[0,T]$, then by (3.37) and the Helly theorem (see [2],p.58,
Th.3.5), we know that there exists a function $\omega\in
BV([0,T];H)$, and a sequence convergent to 0, again denoted by
${\lambda}$, such that
\begin{equation}
\omega_\lambda(t)\rightarrow \omega(t) \  \mathrm{ weakly\  in} \ H
\ \mathrm{for\  every}\  t\in[0,T]
\end{equation}
and
\begin{equation}
\int_t^T\langle\partial\varphi_{\lambda}(y_\lambda(s)), x(s)\rangle
ds\rightarrow\int_t^T\langle d\omega(s),x(s)\rangle,\ \forall x\in
C([t,T];H),\forall t\in[0,T].
\end{equation}
Multiplying equation (3.35) by
sign$p_\lambda(t)=\frac{p_\lambda(t)}{|p_\lambda(t)|}$, we get that
$$\frac{d}{dt}|p_\lambda(t)|=\frac{\nu\|p_\lambda\|^2}{|p_\lambda|}
+\frac{b(p_\lambda,y_\lambda,p_\lambda)}{|p_\lambda|}+ \frac{\langle
\mathscr{C}^*\mathscr{C}(y_\lambda-y^0),p_\lambda\rangle}{|p_\lambda|}
+\frac{\langle\partial\varphi_{\lambda}(y_\lambda),p_\lambda\rangle}{|p_\lambda|}.$$
Since$|b(p_\lambda(t),y_\lambda(t),p_\lambda(t))|\leq
C|p_\lambda(t)|^\frac{1}{2}\|p_\lambda(t)\|^\frac{3}{2}\|y_\lambda(t)\|$,
we get by Young's inequality that
$$\frac{b(p_\lambda(t),y_\lambda(t),p_\lambda(t))}{|p_\lambda(t)|}\leq C\frac{\|p_\lambda(t)\|^\frac{3}{2}}{|p_\lambda(t)|^\frac{3}{4}}|p_\lambda(t)|^\frac{1}{4}\leq\frac{\nu}{2}\frac{\|p_\lambda(t)\|^2}{|p_\lambda(t)|}+C|p_\lambda(t)|.$$
Integrating on $(t,T)$, we obtain by Young's inequality that
$$|p_\lambda(t)|+\frac{\nu}{2}\int_t^T\|p_\lambda(s)\|ds\leq C_1+C_2\int_t^T|p_\lambda(s)|ds.$$
It follows by Gronwall's inequality that
\begin{equation}
\|p_\lambda(t)\|_{L^\infty(0,T;H)}\leq C,
\end{equation}
and so by Alaoglu's theorem,
\begin{equation}
p_\lambda(t)\rightarrow p(t)\ \ \ w^*-L^\infty(0,T;H).
\end{equation}
By (3.32), we infer that
\begin{equation}
\int_t^TU_\lambda(s,t)\mathscr{C}^*\mathscr{C}(y_\lambda(s)-y^0(s))ds\rightarrow
\int_t^TU(s,t)\mathscr{C}^*\mathscr{C}(y^*(s)-y^0(s))ds \
\mathrm{weakly\  in\  }H,
\end{equation}
and by (3.38), we have
\begin{equation}
\int_t^TU_\lambda(s,t)\partial\varphi_{\lambda}(y_\lambda(s))ds\rightarrow
\int_t^TU(s,t)d\omega(s) \ \mathrm{weakly\  in\  H}
\end{equation}
Finally, by (3.41), (3.42) and (3.43), we obtain that
$$p(t)=-\int_t^TU(s,t)(\mathscr{C}^*\mathscr{C}(y^{*}(t)-y^{0}(t)))ds-\int_t^TU(s,t)d\omega(s).$$
This means that $p(t)$ satisfies equation (3.5). Since
$$\int_{0}^{T}\langle\partial\varphi_{\lambda}(y_\lambda(t)), y_\lambda(t)-x(t)\rangle dt
\geq\varphi_{\lambda}(y_\lambda(t))-\varphi_{\lambda}(x(t))\geq
0,\forall x\in \mathcal {K},$$ by (3.39), we can pass $\lambda$ to 0
to get
$$\int_{0}^{T}\langle d\omega(t), y^*(t)-x(t)\rangle \geq 0,$$
so (3.6) holds. To complete the proof, it remains to proof (3.4). By
(3.36) and the definition of $\partial h_\lambda$, we have that
\begin{equation}
\int_0^T(h_\lambda(u_\lambda)-h_\lambda(v))dt\leq \int_0^T\langle
D^*p_\lambda-(u_\lambda-u^*),u_\lambda-v\rangle_Udt, \forall v\in
L^2(0,T;U).
\end{equation}
Remembering that $h_\lambda(u)\leq h(u)$, we obtain by (3.21) and
(3.44) that
\begin{equation}
\int_0^T(h(u^*)-h(v))dt\leq \int_0^T\langle D^*p,u^*-v\rangle_Udt,
\forall v\in L^2(0,T;U).
\end{equation}
This implies the pointwise inequality:
$$\langle
D^*p,u^*-\tilde{v}\rangle_U\geq h(u^*)-h(\tilde{v}),a.e.\
\mathrm{in}\ (0,T).$$ This shows that $D^*p(t)\in \partial
h(u^*(t)),a.e.\ \mathrm{in}\ (0,T)$
$\sharp$ \\
\ \ \\
{\bf Proof of Theorem 3:} Step 1 is the same as the proof of theorem
2. To pass $p_{\lambda}(t),
\partial h_\lambda(u_\lambda(t))$ and $\partial\varphi_\lambda(y_\lambda(t))$ to limit, we need to prove the following
lemma first: \\
{\bf Lemma 3.} \emph{$z_{\lambda}(t) \rightarrow \tilde z(t)\
\mathrm{strongly\ in}\ C([0,T];V)$, where $z_{\lambda}(t)$ and
$\tilde z(t)$ are the solutions to equation (3.22) and equation (3.3) respectively.} \\
{\bf Proof:} We have
\begin{equation}
(z_{\lambda}(t)-\tilde{z}(t))'+\nu
A(z_{\lambda}'(t)-\tilde{z}(t))+(B'(y_{\lambda}(t)))z_{\lambda}(t)-(B'(y^{*}(t)))\tilde{z}(t)
=B(y_{\lambda}(t))-B(y^{*}(t))
\end{equation}
Multiplying equation (3.46) by $A(z_{\lambda}(t)-\tilde{z}(t))$,
integrating on $(0,t)$, we get that
$$\frac{1}{2}\|z_{\lambda}(t)-\tilde{z}(t)\|^2+\nu\int_0^t|A(z_{\lambda}(s)-\tilde{z}(s))|^2ds$$
$$=-\int_0^tb(y_{\lambda}(s)-y^*(s),z_\lambda(s),Az_{\lambda}(s)-A\tilde{z}(s))-\int_0^tb(y^*(s),z_{\lambda}(s)-\tilde{z}(s),Az_{\lambda}(s)-A\tilde{z}(s))
ds$$
$$-b(z_{\lambda}(s)-\tilde{z}(s),y_{\lambda}(s),Az_{\lambda}(s)-A\tilde{z}(s))
-b(\tilde{z}(s),y_{\lambda}(s)-y^*(s),Az_{\lambda}(s)-A\tilde{z}(s))ds$$
$$+\int_0^t\langle B(y_{\lambda}(s))-B(y^{*}(s)),Az_{\lambda}(s)-A\tilde{z}(s)\rangle ds$$
$$\leq
C\int_0^t\left(|y_{\lambda}-y^*|^{\frac{1}{2}}\|Az_\lambda\||A(z_{\lambda}-\tilde{z})|
+\|z_{\lambda}-\tilde{z}\|^{\frac{1}{2}}|Az_{\lambda}-A\tilde{z}|^{\frac{3}{2}}\right)ds$$
$$+C\int_0^t\left(|z_{\lambda}-\tilde{z}|^{\frac{1}{2}}\|z_{\lambda}-\tilde{z}\|^{\frac{1}{2}}|Az_{\lambda}-A\tilde{z}||Ay_\lambda|+\|y_{\lambda}-y^*\|^{\frac{1}{2}}|Ay_{\lambda}-Ay^*|^{\frac{1}{2}}|A(z_{\lambda}-\tilde{z})|\right)ds$$
$$+\int_0^t\langle B(y_{\lambda}(s))-B(y^{*}(s)),Az_{\lambda}(s)-A\tilde{z}(s)\rangle ds$$
$$\leq
C_1\int_0^t\|z_{\lambda}-\tilde{z}\|^2ds+C_2\left(\int_0^t|By_{\lambda}-By^{*}|^2ds
+\|y_{\lambda}-y^*\|_{C([0,T];V)}+\|z_{\lambda}-\tilde{z}\|_{C([0,T;H])}\right)$$
$$+\frac{\nu}{2}\int_0^t|Az_{\lambda}(s)-A\tilde{z}(s))|^2ds$$ Denote
$C_2\left(\int_0^t|B(y_{\lambda}(s))-B(y^{*}(s))|^2ds
+\|y_{\lambda}-y^*\|_{C([0,T];V)}+\|z_{\lambda}-\tilde{z}\|_{C([0,T;H])}\right)$
by $\varepsilon_\lambda$. By the latter inequality and Gronwall's
inequality, it follows that
\begin{equation}
\sup_{t\in[0,T]}\|z_{\lambda}(t)-\tilde{z}(t)\|^2\leq\varepsilon_\lambda
e^{C_1T}\rightarrow 0\ as\ \lambda\rightarrow 0
\end{equation}

Indeed, since $B(y_{\lambda}(t))\rightarrow B(y^{*}(t))$ strongly in
$L^2(0,T;H)$, and $z_{\lambda}\rightarrow\tilde{z}$ strongly in
$C([0,T];H)$, it suffices to prove that $y_{\lambda}\rightarrow y^*$
strongly in $C([0,T];V)$. We have
\begin{equation}
\left\{\begin{array}{lllll}(y_{\lambda}(t)-y^*(t))'+\nu
A(y_{\lambda}(t)-y^*(t))+By_{\lambda}(t)-By^{*}(t)=Du_{\lambda}(t)-Du^{*}(t),\\
y_{\lambda}(0)-y^*(0)=0
\end{array}\right.
\end{equation}
Multiplying equation (3.48) by $Ay_{\lambda}(t)-Ay^*(t)$,
integrating on $(0,t)$. It follows that
$$\frac{1}{2}\|y_{\lambda}(t)-y^*(t)\|^2+\nu\int_0^t|Ay_{\lambda}(s)-Ay^*(s)|^2ds$$
$$\leq
C_1\int_0^t\|y_{\lambda}(t)-y^*(t))\|^2ds+C_2\int_0^t\left(|Du_{\lambda}(s)-Du^{*}(s)|^2+\|y_\lambda-y^*\|_{C([0,T;H])}\right)ds
$$
$$+\frac{\nu}{2}\int_0^t|Ay_{\lambda}(s)-Ay^*(s)|^2ds.$$ Since $u_{\lambda}(t)\rightarrow u^{*}(t)\ \mathrm{strongly\  in}\
L^2(0,T;H)$, $y_\lambda\rightarrow y^*$ strongly in $C([0,T];H)$, we
obtain by Gronwall's inequality that
$$\|y_{\lambda}(t)-y^*(t)\|^2\leq C\left(\|Du_{\lambda}(t)-Du^{*}(t)\|_{L^2(0,T;H)}^2+\|y_\lambda-y^*\|_{C([0,T;H])}\right)dte^{C_1T}
\rightarrow 0\ as\ \lambda\rightarrow 0.$$ This shows that
\begin{equation}
y_{\lambda}\rightarrow y^* \ \mathrm{strongly\  in}\  C([0,T];V) . \
\end{equation}
So (3.47) holds. We complete the proof of Lemma 3.

Applying the similar method, we can get the following result
\begin{equation}
\left\{\begin{array}{lllll}
\|U_\lambda^*(t,s)\|_{L(V,V)}\leq C\\
U_\lambda^*(s,t)\xi\rightarrow U^*(s,t)\xi\ \ \mathrm{in} \
C([0,T];V), \forall \xi\in V
\end{array}\right.
\end{equation}
Now we come back to pass $p_{\lambda}(t),
\partial h_\lambda(u_\lambda(t))$ and $\partial\varphi_\lambda(y_\lambda(t))$ to limit.

By assumption (v) and Lemma 3, we know that $\exists \rho>0,
\lambda_0>0\ s.t.\ z_\lambda(t)+\rho h\in K,\ $ for $t$ in a dense
subset of [0,T], $\forall \lambda<\lambda_0, \forall \|h\|=1$. By
the similar arguments to the proof of Theorem 2, we can get that
\{$\partial\varphi_{\lambda}(y_\lambda)$\} is bounded in
$L^1(0,T;V')$. We denote
$\omega_\lambda(t)=\int_0^t\partial\varphi_{\lambda}(y_\lambda(s))ds,\
t\in[0,T]$, then we know by the Helly theorem that there exist a
functional $\omega\in BV([0,T];V')$, and a sequence convergent to 0,
again denoted by ${\lambda}$, such that
\begin{equation}
\omega_\lambda(t)\rightarrow \omega(t) \ \mathrm{weakly \ in} \ V' \
\mathrm{for \  every}\  t\in[0,T]
\end{equation}
and
\begin{equation}
\int_t^T(\partial\varphi_{\lambda}(y_\lambda(s)), x(s))_{(V',V)}
ds\rightarrow\int_t^T( d\omega(s),x(s))_{(V',V)},\ \forall x\in
C([t,T];V),\forall t\in[0,T].
\end{equation}
Multiply equation (3.35) by $\frac{A^{-1}
p_\lambda(t)}{\|p_\lambda(t)\|_{V'}}$ in the sense of the dual
product between $V'$ and $V$, denote
$q_\lambda(t)=A^{-1}p_\lambda(t)$, we have
$$\frac{d}{dt}\|p_\lambda(t)\|_{V'}=\frac{\nu|p_\lambda(t)|^2}{\|p_\lambda(t)\|_{V'}}
+\frac{b(q_\lambda,y_\lambda,p_\lambda)+b(y_\lambda,q_\lambda,p_\lambda)}{\|p_\lambda(t)\|_{V'}}$$
$$+ \frac{\langle
\mathscr{C}^*\mathscr{C}(y^*(t)-y^0(t)),q_\lambda(t)\rangle}{\|p_\lambda(t)\|_{V'}}
+\frac{\langle\partial\varphi_{\lambda}(y_\lambda),q_\lambda(t)\rangle}{\|p_\lambda(t)\|_{V'}}.$$
Since
$$\frac{b(q_\lambda,y_\lambda,p_\lambda)}{\|p_\lambda\|_{V'}}\leq C\frac{\|q_\lambda\|_{\frac{3}{2}+\varepsilon}|p_\lambda(t)|}{\|p_\lambda(t)\|_{V'}}
\leq
C\|p_\lambda\|_{V'}^{\frac{1-2\varepsilon}{4}}\frac{|p_\lambda|^{\frac{3+2\varepsilon}{2}}}{\|p_\lambda(t)\|_{V'}^{\frac{3+2\varepsilon}{4}}},0<\varepsilon<\frac{1}{2},$$
$$\frac{b(y_\lambda,q_\lambda,p_\lambda)}{\|p_\lambda\|_{V'}}\leq C\|p_\lambda\|_{V'}^{\frac{1}{4}}\frac{|p_\lambda|^{\frac{3}{2}}}{\|p_\lambda(t)\|_{V'}^{\frac{3}{4}}},$$
we obtain by Young's inequality that
$$\|p_\lambda(t)\|_{V'}+\nu\int_t^T\frac{|p_\lambda(s)|^2}{\|p_\lambda(s)\|_{V'}}ds\leq C_1+C_2\int_t^T\|p_\lambda(s)\|_{V'}ds+\frac{\nu}{2}\int_t^T\frac{|p_\lambda(s)|^2}{\|p_\lambda(s)\|_{V'}}ds.$$
It follows by Gronwall's inequality that
\begin{equation}
\|p_\lambda(t)\|_{L^\infty(0,T;V')}\leq C,
\end{equation}
and so by Alaoglu's theorem,
\begin{equation}
p_\lambda(t)\rightarrow p(t)\ \ \ w^*-L^\infty(0,T;V')
\end{equation}
By (3.18) and (3.50), we have
\begin{equation}
\int_t^TU_\lambda(s,t)\mathscr{C}^*\mathscr{C}(y_\lambda(s)-y^0(s))ds\rightarrow
\int_t^TU(s,t)\mathscr{C}^*\mathscr{C}(y(s)-y^0(s))ds \
\mathrm{weakly\  in\  V'.}
\end{equation}
By (3.51), we get
\begin{equation}
\int_t^TU_\lambda(s,t)\partial\varphi_{\lambda}(y_\lambda(s))ds\rightarrow
\int_t^TU(s,t)d\omega(s) \ \mathrm{weakly\  in\ V'.}
\end{equation}
Finally, we obtain by (3.54), (3.55) and (3.56) that
$$p(t)=-\int_t^TU(s,t)(\mathscr{C}^*\mathscr{C}(y^{*}(t)-y^{0}(t)))ds-\int_t^TU(s,t)d\omega(s).$$
This means that $p(t)$ satisfies equation (3.5). Moreover, (3.4) and
(3.8) also hold by applying the same arguments as that in the proof
of theoerm 1. $\sharp$ \\

We shall consider the reflexive Banach space $E$ as $H$ or $V$, and
denote by $(\cdot,\cdot)$ the dual product between $E$ and it's dual
of $E$ (When $E=H$, it is the scalar product in $H$), by $\|\cdot\|$
the norm of $E$. Under the hypothesis of Theorem 2 or the hypothesis
of Theorem 3, We give a corollary
here:\\
{\bf Corollary 1.} \emph{Let the pair $(y^*,u^*)$ be the optimal
pair in problem (P), then there exist $\omega(t)$$\in BV([0,T];E')$
and $p$ satisfying along with $y^*,u^*$, equations (3.4),(3.5),(3.6)
(or(3.8)) and
\begin{equation}
\omega_a(t)\in N_K(y^*(t)), a.e.t\in (0,T)
\end{equation}
\begin{equation}
d\omega_s\in\mathscr{N}_\mathcal {K}(y^*)
\end{equation}}
Here $\omega_a(t)$ is the weak derivative of $\omega(t)$, and
$d\omega_s$ is the singular part of measure $d\omega$. $N_K(y^*(t))$
is the normal cone to $K$ at $y^*(t)$, and $\mathscr{N}_\mathcal
{K}(y^*)$ is the normal cone to $\mathcal {K}$ at $y^*$ which is
precised in definition 1.
\\
{\bf Proof}: Let $t_0$ be arbitrary but fixed in (0,T). For $y\in K$
and $\varepsilon>0$, define the function $y_\varepsilon$
$$
y_\varepsilon(t)=\left\{\begin{array}{lllll}
y^*(t)\ \ \ \ \ \ \ \ \ \ \ \ \mathrm{for}\ |t-t_0|\geq\varepsilon\\
(1-\varepsilon^{-1}(t_0-t))y+\varepsilon^{-1}(t_0-t)y^*(t_0-\varepsilon)
\ \ \ \ \ \mathrm{for}\ \ t\in[t_0-\varepsilon,t_0]\\
(1-\varepsilon^{-1}(t-t_0))y+\varepsilon^{-1}(t-t_0)y^*(t_0+\varepsilon)
\ \ \ \ \ \mathrm{for}\ \ t\in[t_0,t_0+\varepsilon].\\
\end{array}
\right. $$ Obviously $y_\varepsilon$ is continuous from [0,T] to $E$
and $y_\varepsilon(t)\in K, \forall t\in [0,T]$. By (4.6)(or (4.8)),
we have
\begin{equation}
\int_0^T(\dot{\omega}(t),y^*(t)-y_\varepsilon(t))dt+\int_0^T(d\omega_s,y^*-y_\varepsilon)\geq0.
\end{equation}
We set
$\rho_\varepsilon(t)=\varepsilon^{-1}(y^*(t)-y_\varepsilon(t))$. If
$t_0$ happens to be a Lebesgue point for the function $\omega_a$,
then by an elementary calculation involving the definition of
$y_\varepsilon$, we get
\begin{equation}
\lim_{\varepsilon\rightarrow
0}\int_0^T(\omega_a(t),\rho_\varepsilon(t))dt=(\omega_a(t_0),y^*(t_0)-y).
\end{equation}
Inasmuch as $y^*(t)-y_\varepsilon(t)=0$ outside $[t_0-\varepsilon,
t_0+\varepsilon]$, we have
$$\int_0^T(d\omega_s,y^*-y_\varepsilon)=\int_{t_0-\varepsilon}^{t_0+\varepsilon}(d\omega_s,y^*-y_\varepsilon).$$
On the other hand, for each $\eta >0$, there exists
$\{y^*_{i\eta}\}_{i=1}^{N}\subset E$, and $\alpha_{i\eta}\in
C([0,T])$ such that
$$\|y^*(t)-\sum_{i=1}^Ny^*_{i\eta}\alpha_{i\eta}(t)\|\leq\eta\ \mathrm{for}\ t\in [0,T]$$
We set $z_\eta (t)=y^*(t)-\sum_{i=1}^Ny^*_{i\eta}\alpha_{i\eta}(t)$,
then we have
$$|\int_{t_0-\varepsilon}^{t_0+\varepsilon}(d\omega_s,z_\eta)|\leq(V_s(t_0+\varepsilon)-V_s(t_0-\varepsilon))\cdot \sup\{\|z_\eta\|;|t-t_0|\leq \varepsilon\}$$
where $V_s(t)$ is the variation of $\omega_s$ on the interval [0,t].
Since $V_s$ is a.e. differentiable on (0,T), we may assume that
\begin{equation}
\limsup_{\varepsilon\rightarrow
0}\varepsilon^{-1}\int_{t_0-\varepsilon}^{t_0+\varepsilon}(d\omega_s,z_\eta)\leq
C\eta,
\end{equation}
where $C$ is independent of $\eta$. Now, we have
$$|\int_{t_0-\varepsilon}^{t_0+\varepsilon}(d\omega_s,\sum_{i=1}^Ny^*_{i\eta}\alpha_{i\eta})|\leq\sum_{i=1}^N|\int_{t_0-\varepsilon}^{t_0+\varepsilon}\alpha_{i\eta}d(\omega_s,y^*_{i\eta})|$$
$$\leq\sum_{i=1}^N(V_{i\eta}(t_0+\varepsilon)-V_{i\eta}(t_0-\varepsilon))\gamma_{i\eta},$$
where $V_{i\eta}(t)$ is the variation of $(\omega_s,y^*_{i\eta})$ on
interval [0,t] and $\gamma_{i\eta}=\sup|\alpha_{i\eta}(t)|$. Since
the weak derivative of $\omega_s$ is zero a.e. on (0,T). We may
infer that:
$$\frac{d}{dt}V_{i\eta}(t)=0,\ a.e.
 \ \mathrm{in}\ (0,T),$$
and therefore we may assume that
\begin{equation}
\lim_{\varepsilon\rightarrow
0}\varepsilon^{-1}\int_{t_0-\varepsilon}^{t_0+\varepsilon}(d\omega_s,\sum_{i=1}^Ny^*_{i\eta}\alpha_{i\eta}(t))=0,\
\forall \eta>0.
\end{equation}
By (3.61) and (3.62) we have
$$\lim_{\varepsilon\rightarrow
0}\varepsilon^{-1}\int_0^T(d\omega_s,\rho_\varepsilon)=0.$$ By
(3.59) and (3.60) we get
\begin{equation}
(\omega_a(t_0),y^*(t_0)-y)\geq0, \ \ a.e.\ t_0\in(0,T).
\end{equation}
Since $y$ is arbitrary, (3.57) holds.

To conclude the proof it remains to be shown that
$d\omega_s\in\mathscr{N}_\mathscr{K}(y^*)$, that is
$\int_0^T(d\omega_s,y^*-x)\geq0,\ \forall x\in \mathscr{K}$. Let
$\mathcal {O}$ be the support of the singular measure $d\omega_s$.
Then for any $\varepsilon>0$, there exists an open subset $\mathcal
{U}$ of (0,T), s.t. $\mathcal {O}\subset\mathcal {U}$ and
$m(\mathcal {U})\leq\varepsilon$, where $m$ is the lebesgue measure.
Let $\rho\in C_0^\infty(\mathbb{R})$ be such that
$0\leq\rho\leq1,\rho=1$ on $\mathcal {O}$ and $\rho=0$ on
$(0,T)\setminus\mathcal {U}$. We set $y^\varepsilon=\rho
x+(1-\rho)y^*$, where $x\in \mathscr{K}$ is arbitrary. By
(3.6)(or(3.8)), we have
\begin{equation}
\int_0^T(\omega_a(t),y^*(t)-y^\varepsilon(t))dt+\int_0^T(d\omega_s,y^*-y^\varepsilon)\geq0.
\end{equation}
Since $y^*-y^\varepsilon=0$ on $(0,T)\setminus\mathcal {U}$, we
obtain that
$$|\int_0^T(\omega_a(t),y^*(t)-y^\varepsilon)dt|\leq \int_\mathcal {U}|\omega_a|dt\leq\delta(\varepsilon),$$
where $\lim_{\varepsilon\rightarrow 0}\delta(\varepsilon)=0$. On the
other hand, since $d\omega_s=0$ on $(0,T)\setminus\mathcal {O}$ and
$\rho=1$ on $\mathcal {O}$, we see that
$$\int_0^T(d\omega_s,y^*-y^\varepsilon)=\int_0^T(d\omega_s,y^*-x).$$
It follows that $\int_0^T(d\omega_s,y^*-x)\geq
-\delta(\varepsilon)$. Since $\varepsilon$ is arbitrary, (3.58)
holds. $\sharp$

\section{Examples}\setcounter{equation}{0}
In this section, we shall give some applications of the above
results in some special cases of state constraints wherein Theorem 1 and Theorem 2 can be applied.\\
{\bf Example 1.} Let $K$ be the set $K=\{y\in H;
\int_\Omega|y(x)|^2dx\leq \rho^2\}$, then $K$ is a closed convex set
in $H$, since
$$\|\tilde{z}(t)\|_{C([0,T];H)}\leq C(\|B(y^*(t))+D\tilde{u}(t)+f(t)\|_{L^2(0,T;H))})$$
so it is feasible to apply theorem 2 to get the necessary condition
of the optimal control pair after checking whether condition (iv) is
satisfied or not. The set $K$ physically gives a constraint on the
turbulence kinetic energy, which is usually bounded instead of
infinite. In this case, the maximum principle can be described as
following:
\begin{equation}
D^*p(t)\in\partial h(u^*(t)) \ \ a.e. [0,T]
\end{equation}
\begin{equation}
\noindent\left\{\begin{array}{lllll}
     p_1'(t)=\nu Ap_1(t)+(B'(y^{*}(t))^{*})p_1(t)+\mathscr{C}^*\mathscr{C}(y^{*}(t)-y^{0}(t))+\omega_a(t),\ a.e.\ \mathrm{in}\ (0,T)\\
     p_1(T)=0
\end{array}\right.
\end{equation}
\begin{equation}
\noindent\left\{\begin{array}{lllll}
     p'_2(t)=\nu Ap_2(t)+(B'(y^{*}(t))^{*})p_2(t)+d\omega_s,\ a.e.\ \mathrm{in}\ (0,T)\\
     p_2(T)=0
\end{array}\right.
\end{equation}
Moreover,
\begin{equation}
\omega_a(t)\in N_K(y^*(t))=\{\lambda(t)y^*(t);\lambda(t)\geq 0,
a.e.\ \mathrm{in}\ (0,T)\}
\end{equation}
Here $p_1,p_2$ is the decomposition of $p$, that is
$p(t)=p_1(t)+p_2(t)$. Since $\omega_a(t)\in L^1(0,T;H)$, $d\omega_s$
is the singular part of the measure $d\omega$, we know that equation
(4.2) has a strong solution $p_1\in C([0,T];H)$, while equation
(4.3) has only a mild solution $p(t)=-\int_t^TU(s,t)dw_s$.\\
{\bf Example 2.} Let $K$ be the so called Enstrophy set
$$K=\{y\in V; \int_\Omega|\nabla\times y|^2dx\leq\rho^2\}$$
where $\nabla\times y=curl\ y(x)$, and it is true that the norm
$|\nabla\times y|$ is equivalent to the norm $\|y\|$ in the space
$V$. In fluid mechanics, the enstrophy $\mathcal
{E}(y)=\int_\Omega|\nabla\times y|^2dx$ can be interpreted as
another type of potential density. More precisely, the quantity
directly related to the kinetic energy in the flow model that
corresponds to dissipation effects in the fluid. It is particularly
useful in the study of turbulent flows, and is often identified in
the study of trusters, as well as the flame field. Enstrophy set
gives a constraint on the vorcity of the fluid motion. Since
$$\|\tilde{z}(t)\|_{C([0,T];V)}\leq C(\|B(y^*(t))+D\tilde{u}(t)+f(t)\|_{L^2(0,T;H))})$$
it is feasible to apply theorem 3 to get the necessary condition of
the optimal control pair after checking whether condition (v) is
satisfied or not. In this case, the maximum principle can be
described by (4.1), (4.2) and (4.3). Moreover,
$$\omega_a(t)\in N_K(y^*(t))=\{\lambda(t)Ay^*(t);\lambda(t)\geq 0, a.e.\ \mathrm{in}\ (0,T)\}$$
{\bf Example 3.} Let $K$ be the so called Helicity set,
$$K=\{y\in V; \int_\Omega\langle y,curl \ y\rangle^2dx+\lambda\int_\Omega|\nabla y|^2dx\leq\rho^2\}$$
where $\lambda$, $\rho$ are positive constants. In fluid mechanics,
helicity is the extent to which corkscrew-like motion occurs. If a
parcel of fluid is moving, undergoing solid body motion rotating
about an axis parallel to the direction of motion, it will have
helicity. If the rotation is clockwise when viewed from ahead of the
body, the helicity will be positive, if counterclockwise, it will be
negative. Helicity is a useful concept in theoretical descriptions
of turbulence. Formally, helicity is defined as
$$H=\int_\Omega\langle y,curl \ y\rangle dx$$
The helicity set plays an important role in fluid mechanics, and in
particular, it is an invariant set of Euler's equation for
incompressible fluids(See [4]). This set gives a constraint on the
helicity and the smoothness of the velocity field. By the same
argument as in Example 2, we know that it is feasible to apply
theorem 3 to get the necessary condition of the optimal pair when
the state constrained set is Helicity set, and in this case, the
maximum principle can be described by (4.1), (4.2), (4.3). Moreover,
$$\omega_a(t)\in N_K(y^*(t))=\{\lambda(t)(Ay^*(t)+\mathrm{curl}y^*);\lambda(t)\geq 0,a.e.\ \mathrm{in}\ (0,T)\}$$\\
\ \ \ \\

\bigskip

\begin{thebibliography}{00}
\bibitem{01}\label{Barbu1}V.Barbu, Analysis and Control of Nonlinear Infinite Dimensional System, Academic Press, Boston, 1993.
\bibitem{01}\label{Barbu2}V.Barbu, T.Precupanu, Convexity and Optimization in Banach Spaces,
Math.Appl.(East European Series), Vol.10, Romanian edition,
D.Reidel, Dordrecht(1986).
\bibitem{00}\label{Barbu3}V.Barbu, Nonlinear Differential Equations of Monotone Type in Banach Spaces, Springer-Verlag, Berlin, 2010.
\bibitem{00}\label{Barbu4}V.Barbu, N.H.Pavel, Flow-invariant closed
sets with respect to nonlinear semigroup flows, Nonlinear differ.
equ. appl. 10(2003) 57-72.
\bibitem{00}\label{Barbu5}V.Barbu, Optimal control of Navier-Stokes equations with periodic inputs. Nonlinear Anal. TMA 31(1/2)(1998)15-31.
\bibitem{00}\label{Thomas1}Thomas R.Bewley, Flow Control: new challenges for a new
Renaissance, Progress in Aerospace Sciences 37(2001)21-58.
\bibitem{00}\label{Thomas2}TR.Bewley, R.Temam, M.Ziane, Existence
and uniqueness of optimal control to the Navier-Stokes equations,
C.R.Acad.Sci.Paris,t.330,Serie 1,p.1-5,2000.
\bibitem{00}\label{Thomas3}TR.Bewley, R.Temam, M.Ziane, A general
framework for robust control in fluid mechanics, Physica D 138(2000)
360-392.
\bibitem{00}\label{Temman}R.Temman, Navier-Stokes Equations, North-Holland, Amsterdam,
1979.
\bibitem{00}\label{Wang}G.Wang, L.Wang, Maximum principle of state-constrained optimal
control governed by fluid dynamic systems, Nonlinear Analysis,
52(2003) 1911-1931.
\bibitem{00}\label{Wang}G.Wang, Optimal controls of 3-dimensional Navier-Stokes equations with state constraints.
 SIAM.J.CONTROL OPTIM.Vol.41(2002), No.2,pp.583-606.
\end{thebibliography}
\end{document}